\documentclass[12pt]{amsart}      
\usepackage{amssymb}
\usepackage{eucal}
\usepackage{amsmath}
\usepackage{amscd}
\usepackage[all]{xy}           

\usepackage{amsfonts,latexsym}
\usepackage{xspace}
\usepackage{hyperref}
\usepackage{float}
\usepackage{color}
\usepackage{colordvi}
\usepackage{multicol}

\newcommand{\nc}{\newcommand}
\newcommand{\delete}[1]{}

\nc{\dfootnote}[1]{{}}          
\nc{\ffootnote}[1]{\dfootnote{#1}}

\delete{
\nc{\mfootnote}[1]{{}}        
\nc{\ofootnote}[1]{{}}        
}

\nc{\mfootnote}[1]{\footnote{#1}} 
\nc{\ofootnote}[1]{\footnote{\tiny Older version: #1}} 

\nc{\mlabel}[1]{\label{#1}}  
\nc{\mcite}[1]{\cite{#1}}  
\nc{\mref}[1]{\ref{#1}}  
\nc{\mkeep}[1]{{}}      
\nc{\mbibitem}[1]{\bibitem{#1}} 

\delete{
\nc{\mcite}[1]{\cite{#1}{{\bf{{\ }(#1)}}}}  
\nc{\mlabel}[1]{\label{#1}  
{\hfill \hspace{1cm}{\bf{{\ }\hfill(#1)}}}}
\nc{\mref}[1]{\ref{#1}{{\bf{{\ }(#1)}}}}  
\nc{\mbibitem}[1]{\bibitem[\bf #1]{#1}} 
\nc{\mkeep}[1]{\marginpar{{\bf #1}}} 
}

\newtheorem{theorem}{Theorem}[section]
\newtheorem{prop}[theorem]{Proposition}

\newtheorem{lemma}[theorem]{Lemma}

\newtheorem{prop-def}{Proposition-Definition}[section]

\nc{\bond}{\vdash}
\nc{\mtiny}{\scriptscriptstyle}
\nc{\comp}[1]{\langle #1\rangle} \nc{\spr}{\cdot}
\nc{\disp}[1]{\displaystyle{#1}}
\nc{\sk}[1]{\mrm{sk}(#1)}
\nc{\ve}[1]{\mrm{vec}({#1})}
\nc{\pf}[1]{\check{#1}}
\nc{\Li}{\mrm{Li}}
\nc{\icut}{^!}
\nc{\bin}[2]{ (_{\stackrel{\scs{#1}}{\scs{#2}}})}  
\nc{\binc}[2]{ \bigg (\!\! \begin{array}{c} \scs{#1}\\
    \scs{#2} \end{array}\!\! \bigg )}  
\nc{\bbinc}[2]{ \left (\!\! \begin{array}{c} {#1}\\
    {#2} \end{array}\!\! \right )}  
\nc{\bincc}[2]{  \left ( {\scs{#1} \atop
    \vspace{-.5cm}\scs{#2}} \right )}  
\nc{\sarray}[2]{\begin{array}{c}#1 \vspace{.1cm}\\ \hline
    \vspace{-.35cm} \\ #2 \end{array}}
\nc{\spair}[2]{\big[\begin{array}{c}\scs{#1} \\ \scs{#2} \end{array} \big]}
\nc{\ppair}[2]{\big\langle\begin{array}{c}\scs{#1} \\ \scs{#2} \end{array} \big\rangle}
\nc{\zsg}[1]{\widehat{#1}} \nc{\bs}{\bar{S}}
\nc{\dcup}{\stackrel{\bullet}{\cup}}
\nc{\dbigcup}{\stackrel{\bullet}{\bigcup}} \nc{\la}{\longrightarrow}
\nc{\fe}{\'{e}} \nc{\rar}{\rightarrow} \nc{\dar}{\downarrow}
\nc{\dap}[1]{\downarrow \rlap{$\scriptstyle{#1}$}}
\nc{\uap}[1]{\uparrow \rlap{$\scriptstyle{#1}$}}
\nc{\dt}[1]{{#1}^\sharp} \nc{\st}[1]{{#1}^\flat}
\nc{\defeq}{\stackrel{\rm def}{=}} \nc{\dis}[1]{\displaystyle{#1}}
\nc{\dotcup}{\ \displaystyle{\bigcup^\bullet}\ } \nc{\hcm}{\
\hat{,}\ } \nc{\hcirc}{\hat{\circ}} \nc{\hts}{\hat{\shpr}}
\nc{\lts}{\stackrel{\leftarrow}{\shpr}}
\nc{\rts}{\stackrel{\rightarrow}{\shpr}} \nc{\lleft}{[}
\nc{\lright}{]} \nc{\uni}[1]{\tilde{#1}} \nc{\free}[1]{\bar{#1}}
\nc{\den}[1]{\check{#1}} \nc{\lrpa}{\wr} \nc{\curlyl}{\left \{
\begin{array}{c} {} \\ {} \end{array}
    \right . \!\!\!\!\!\!\!}
\nc{\curlyr}{ \!\!\!\!\!\!\!
    \left . \begin{array}{c} {} \\ {} \end{array}
    \right \} }
\nc{\longmid}{\left | \begin{array}{c} {} \\ {} \end{array}
    \right . \!\!\!\!\!\!\!}
\nc{\ot}{\otimes} \nc{\bigot}{\bigotimes} \nc{\mdiv}{\mrm{div}}
\nc{\shpf}{\frakp} \nc{\sg}{G} \nc{\ssg}[1]{\overline{#1}}
\nc{\psg}[1]{\widetilde{#1}} \nc{\zg}{{Z}} \nc{\ig}{I} \nc{\pg}{P}
\nc{\jg}{J} \nc{\eg}{E} \nc{\fg}{F} \nc{\cg}{C} \nc{\mg}{M}
\nc{\abg}{C} \nc{\ug}{{U}} \nc{\bas}{B} \nc{\Lyn}{\mrm{Lyn}}
\nc{\lyn}{\Lyn} \nc{\sG}{{\cals}} \nc{\iG}{{\cali}}
\nc{\jG}{{\calj}} \nc{\eG}{{\cale}} \nc{\pG}{{\calp}}
\nc{\fG}{{\calf}} \nc{\cG}{{\calc}} \nc{\mG}{{\calm}}
\nc{\ora}[1]{\stackrel{#1}{\rar}}
\nc{\ola}[1]{\stackrel{#1}{\la}}
\nc{\pex}[1]{\{#1\}} \nc{\scs}[1]{\scriptstyle{#1}}
\nc{\mrm}[1]{{\rm #1}} \nc{\sym}[1]{{\widehat{#1}}}
\nc{\margin}[1]{\marginpar{\rm #1}}   
\nc{\dirlim}{\displaystyle{\lim_{\longrightarrow}}\,}
\nc{\invlim}{\displaystyle{\lim_{\longleftarrow}}\,}
\nc{\mvp}{\vspace{0.5cm}} \nc{\svp}{\vspace{2cm}}
\nc{\vp}{\vspace{8cm}} \nc{\proofbegin}{\noindent{\bf Proof: }}
\nc{\proofend}{$\blacksquare$ \vspace{0.5cm}}
\nc{\shqs}{\eta} \font\cyr=wncyr10
\nc{\sha}{{\mbox{\cyr X}}}  
\newfont{\scyr}{wncyr10 scaled 550}
\nc{\ssha}{\,\mbox{\bf \scyr X}\,}
\newfont{\bcyr}{wncyr10 scaled 1000}
\nc{\qssha}{{{\ssha\hspace{-2pt}_\ast}}\,} \nc{\pssha}{{\star}}
\nc{\qsshab}{{{\ssha\hspace{-2pt}_{\rho}}}\,}
\nc{\ncsha}{{\mbox{\cyr X}^{\mathrm NC}}} \nc{\ncshao}{{\mbox{\cyr
X}^{\mathrm NC,\,0}}}
\nc{\shpr}{\diamond}    

\nc{\shf}{{^{\ssha}}} \nc{\qsh}{{^{\ast}}}
\nc{\psh}{{^{\pi}}}
\nc{\pfsha}{{\ssha_{\pi}}}
\nc{\esh}{{^{ext}}}
\nc{\lshf}{_{\ssha}} \nc{\lqsh}{_{\ast}} \nc{\lzero}{_{\hskip -5pt
0}} \nc{\shzero}{_{\hskip -7.5pt 0}} \nc{\lone}{_{\hskip -7.5pt 1}}
\nc{\shprl}{{{\shpr}_\lambda}}
\nc{\shpro}{\diamond^0}    
\nc{\shpru}{\check{\diamond}} \nc{\catpr}{\diamond_l}
\nc{\rcatpr}{\diamond_r} \nc{\lapr}{\diamond_a}
\nc{\lepr}{\diamond_e} \nc{\tcon}{^{\ot}} \nc{\conv}{_c}
\nc{\vep}{\varepsilon} \nc{\labs}{\mid\!} \nc{\rabs}{\!\mid}
\nc{\hsha}{\widehat{\sha}} \nc{\lsha}{\stackrel{\leftarrow}{\sha}}
\nc{\rsha}{\stackrel{\rightarrow}{\sha}}
\nc{\EDS}{{\mrm{EDS}}\xspace} \nc{\DS}{{\mathbf{DS}}}
\nc{\lc}{[} \nc{\rc}{]} \nc{\rbset}{R} \nc{\rbnum}{r}
\nc{\rbfun}{\mathbf{R}} \nc{\pset}{P} \nc{\pnum}{p}
\nc{\pfun}{\mathbf{P}} \nc{\spset}{SP} \nc{\spnum}{sp}
\nc{\spgen}{\mathbf{SP}} \nc{\srbi}[1]{\{#1\}}

\nc{\sumx}{{\fraks_1}} \nc{\sumy}{{\fraks_2}} \nc{\sumz}{{\fraks_3}}

\nc{\sume}{{\fraks_1}} \nc{\suma}{{\fraks_{2,1}}} \nc{\sumb}{{\fraks_{2,2}}}
\nc{\sumc}{{\fraks_{3,1}}} \nc{\sumd}{{\fraks_{3,2}}} \nc{\sumh}{{\fraks_1}}
\nc{\sumi}{{\fraks_2}} \nc{\sumj}{{\fraks_3}}

\nc{\rind}{r} \nc{\sind}{s} \nc{\tind}{t}  \nc{\kdim}{k}
\nc{\ldim}{\ell}

\nc{\alga}{{A}} \nc{\ann}{\mrm{ann}} \nc{\Aut}{\mrm{Aut}}
\nc{\can}{\mrm{can}} \nc{\colim}{\mrm{colim}} \nc{\Cont}{\mrm{Cont}}
\nc{\rchar}{\mrm{char}} \nc{\cok}{\mrm{coker}} \nc{\dtf}{{R-{\rm
tf}}} \nc{\dtor}{{R-{\rm tor}}}

\nc{\Div}{{\mrm Div}} \nc{\End}{\mrm{End}} \nc{\Ext}{\mrm{Ext}}
\nc{\Fil}{\mrm{Fil}} \nc{\Frob}{\mrm{Frob}} \nc{\Gal}{\mrm{Gal}}
\nc{\GL}{\mrm{GL}} \nc{\lord}{\mrm{L-order}\xspace}
\nc{\rme}{\mrm{E}} \nc{\rmt}{\mrm{T}} \nc{\Sym}{\mrm{Sym}}
\nc{\Hom}{\mrm{Hom}} \nc{\hsr}{\mrm{H}} \nc{\hpol}{\mrm{HP}}
\nc{\id}{\mrm{id}} \nc{\im}{\mrm{im}} \nc{\incl}{\mrm{incl}}
\nc{\length}{\mrm{length}} \nc{\leng}{\mrm{\ell}} \nc{\LR}{\mrm{LR}}
\nc{\mchar}{\mrm char}
\nc{\MZV}{\mrm{MZV}\xspace}
\nc{\MZVs}{\mrm{MZVs}\xspace}
\nc{\MZVf}{\mrm{MZV fractions}\xspace}
\nc{\PF}{\mathbf{PF}}
\nc{\MPV}{\mrm{MPV}\xspace}
\nc{\MPVs}{\mrm{MPVs}\xspace}
\nc{\MPL}{\mrm{MPL}\xspace}
\nc{\MPLs}{\mrm{MPLs}\xspace} \nc{\mzvalg}{\mathbf{MZV}}
\nc{\mplalg}{\mathbf{MPV}}
\nc{\pfalg}{\mathbf{PF}}
\nc{\edsalg}{\mathbf{EDS}} \nc{\qeds}{$\QQ$-EDS\xspace}
\nc{\zeds}{$\ZZ$-EDS\xspace} \nc{\zpeds}{$\ZZ_p$-EDS\xspace}
\nc{\fpeds}{$\FF_p$-EDS\xspace} \nc{\NC}{\mrm{NC}}
\nc{\mpart}{\mrm{part}} \nc{\os}{\mrm{OS}} \nc{\qs}{\mrm{QS}}
\nc{\ql}{{\QQ_\ell}} \nc{\qp}{{\QQ_p}} \nc{\rank}{\mrm{rank}}
\nc{\rcot}{\mrm{cot}} \nc{\rdef}{\mrm{def}} \nc{\rdiv}{{\rm div}}
\nc{\rtf}{{\rm tf}} \nc{\rtor}{{\rm tor}} \nc{\res}{\mrm{res}}
\nc{\sh}{\mrm{Sh}} \nc{\TL}{\mrm{TL}} \nc{\Spec}{\mrm{Spec}}
\nc{\tor}{\mrm{tor}} \nc{\Tr}{\mrm{Tr}} \nc{\tr}{\mrm{tr}}
\nc{\ETC}{\mathrm{ETC}} \nc{\ETL}{\mathrm{ETL}}
\nc{\EL}{\mathrm{EL}} \nc{\RETL}{\mathrm{RETL}}
\nc{\EETL}{\widetilde{\TL}} \nc{\word}{\rm word\xspace}
\nc{\words}{\rm words\xspace} \nc{\varab}{\phi_{\alpha,\beta}}

\nc{\ab}{\mathbf{Ab}} \nc{\Alg}{\mathbf{Alg}}
\nc{\Algo}{\mathbf{Alg}^0} \nc{\Bax}{\mathbf{Bax}}
\nc{\Baxo}{\mathbf{Bax}^0} \nc{\RBo}{\mathbf{RB}^0}
\nc{\BRB}{\mathbf{RB}} \nc{\Dend}{\mathbf{DD}} \nc{\bfe}{{\bf e}}
\nc{\bff}{{\bf f}} \nc{\bfk}{{\bf k}} \nc{\bfone}{{\bf 1}}
\nc{\base}[1]{{a_{#1}}} \nc{\detail}{\marginpar{\bf More detail}
    \noindent{\bf Need more detail!}
    \svp}
\nc{\Diff}{\mathbf{Diff}} \nc{\gap}{\marginpar{\bf
Incomplete}\noindent{\bf Incomplete!!}
    \svp}
\nc{\FMod}{\mathbf{FMod}} \nc{\RB}{\mathbf{RB}}
\nc{\Int}{\mathbf{Int}} \nc{\Mon}{\mathbf{Mon}}
\nc{\remarks}{\noindent{\bf Remarks: }} \nc{\Rep}{\mathbf{Rep}}
\nc{\Rings}{\mathbf{Rings}} \nc{\Sets}{\mathbf{Sets}}
\nc{\DT}{\mathbf{DT}}

\nc{\BA}{{\Bbb A}} \nc{\CC}{{\Bbb C}} \nc{\DD}{{\Bbb D}}
\nc{\EE}{{\Bbb E}} \nc{\FF}{{\Bbb F}} \nc{\GG}{{\Bbb G}}
\nc{\HH}{{\Bbb H}} \nc{\LL}{{\Bbb L}} \nc{\NN}{{\Bbb N}}
\nc{\QQ}{{\Bbb Q}} \nc{\RR}{{\Bbb R}} \nc{\TT}{{\Bbb T}}
\nc{\VV}{{\Bbb V}} \nc{\ZZ}{{\Bbb Z}}

\nc{\cala}{{\mathcal A}} \nc{\calc}{{\mathcal C}}
\nc{\cald}{{\mathcal D}} \nc{\cale}{{\mathcal E}}
\nc{\calf}{{\mathcal F}} \nc{\calg}{{\mathcal G}}
\nc{\calh}{{\mathcal H}} \nc{\cali}{{\mathcal I}}
\nc{\calj}{{\mathcal J}} \nc{\call}{{\mathcal L}}
\nc{\calm}{{\mathcal M}} \nc{\caln}{{\mathcal N}}
\nc{\calo}{{\mathcal O}} \nc{\calp}{{\mathcal P}}
\nc{\calr}{{\mathcal R}} \nc{\calt}{{\mathcal T}}
\nc{\calw}{{\mathcal W}} \nc{\calx}{{\mathcal X}}
\nc{\CA}{\mathcal{A}}

\nc\indI{\mathcal{I}}

\nc{\fraka}{{\mathfrak a}} \nc{\frakB}{{\mathfrak B}}
\nc{\frakb}{{\mathfrak b}} \nc{\frakd}{{\mathfrak d}}
\nc{\frakF}{{\mathfrak F}} \nc{\frakf}{{\mathfrak f}}
\nc{\frakg}{{\mathfrak g}} \nc{\frakL}{{\mathfrak L}}
\nc{\frakm}{{\mathfrak m}} \nc{\frakM}{{\mathfrak M}}
\nc{\frakMo}{{\mathfrak M}^0} \nc{\frakp}{{\mathfrak p}}
\nc{\fraks}{{\mathfrak S}}
\nc{\frakw}{{\mathfrak w}}
\nc{\frakx}{{\mathfrak x}} \nc{\ox}{\overline{\frakx}}
\nc{\frakX}{{\mathfrak X}} \nc{\fraky}{{\mathfrak y}}

\nc{\li}[1]{\textcolor{blue}{Li: #1}}
\nc{\wys}[1]{\textcolor{green}{William: #1}}
\nc{\byhs}[1]{\textcolor{red}{Bingyong: #1}}

\nc{\rrb}[1]{[#1]} \nc{\rrrb}[1]{\{#1\}} \nc{\ideal}[1]{\langle
#1\rangle} \nc{\refl}[1]{\overline{#1}} \nc{\rrB}{{reflexive
Rota-Baxter}\xspace} \nc{\rrrB}{{radical reflexive
Rota-Baxter}\xspace} \nc{\srb}{{strict Rota-Baxter}\xspace}
\nc{\Srb}{{Strict Rota-Baxter}\xspace}

\renewcommand\geq{\geqslant}
\renewcommand\leq{\leqslant}

\nc\rbop{{\lc\,\,\rc}} \nc\rbopi[1]{{\lc_{#1} \, \rc_{#1}}}

\nc{\redtext}[1]{\textcolor{red}{#1}}

\begin{document}

\title[Weighted sum formula for MZVs]
{Weighted sum formula for multiple zeta values}
%
\author{Li Guo}
\address{Department of Mathematics and Computer Science,
         Rutgers University,
         Newark, NJ 07102, USA}
\email{liguo@newark.rutgers.edu}
\author{Bingyong Xie}
\address{Department of Mathematics, Peking University, Beijing, 100871, China}
\email{byhsie@math.pku.edu.cn}


\begin{abstract}
The sum formula is a basic identity of multiple zeta values that expresses a Riemann zeta value as a homogeneous sum of multiple zeta values of a given dimension. This formula was already known to Euler in the dimension two case, conjectured in the early 1990s for higher dimensions and then proved by Granville and Zagier independently. Recently a weighted form of Euler's formula was obtained by Ohno and Zudilin. We generalize it to a weighted sum formula for multiple zeta values of all dimensions.

\delete{
\noindent
MSC classes: 11M41, 11M99, 40B05.
\smallskip

\noindent
Keywords: Multiple zeta values, sum formula, Euler's sum formula, weighted sum formula, double shuffle relation.
}
\end{abstract}

\maketitle




\setcounter{section}{0}


\section{Introduction}
Multiple zeta values (\MZVs) are special values of the multi-variable analytic function
\begin{equation}
\zeta(s_1,\cdots, s_k)=\sum_{n_1>\cdots>n_k>0} \frac{1}{n_1^{s_1}\cdots n_k^{s_k}}
\mlabel{eq:mzvs}
\end{equation}
at integers $s_1\geq 2, s_i\geq 1, 1\leq i\leq k.$ Their
study in the two variable case went back to Goldbach and Euler~\mcite{JW}.
The general concept was introduced in the early 1990s, with motivation from both mathematics~\mcite{Ho0,Za} and physics ~\mcite{BK}. Since then
the subject has been studied intensively with interactions to a broad range of areas in mathematics and physics, including arithmetic geometry, combinatorics, number theory, knot theory, Hopf algebra, quantum field theory and mirror symmetry~\mcite{3BL,BB,GM,GZ,Ho1,Ho3,IKZ,LM,Ra,Te}.

A principle goal in the theoretical study of \MZVs is to determine all
algebraic relations among them. Conjecturally all such relations
come from the so-called extended double shuffle relations.
But there are no definite way to exhaust all of them and
new identities of \MZVs are being found steadily~\mcite{BB,ELO,GX2,HO,OZ,Zh3}.

One of the first established and most well-known among these identities is the striking {\bf sum formula}, stating that, for given positive integers $k$ and $n\geq k+1$,
\begin{equation}
\sum_{s_i\geq 1, s_1\geq 2, s_1+\cdots +s_k=n} \zeta (s_1,\cdots,s_k)=\zeta(n),
\mlabel{eq:sum}
\end{equation}
suggesting intriguing connection between Riemann zeta values and multiple zeta values.
This formula was obtained by Euler when $k=2$, known as {\bf Euler's sum formula}~\mcite{Eu}
\begin{equation}
\sum_{i=2}^{n-1} \zeta(i,n-i)=\zeta(n)
\mlabel{eq:euler}
\end{equation}
which includes the basic example $\zeta(2,1)=\zeta(3)$. Its general
form was conjectured in~\mcite{Ho1} and proved by
Granville~\mcite{Gr} and Zagier~\mcite{Za2}. Since then the sum
formula has been generalized and extended in various
directions~\mcite{Br,ELO,HO,KTY,Oh,OO,OW,OZa,OU}.
Ohno and
Zudilin~\mcite{OZ} recently proved a weighted form of Euler's sum
formula (the {\bf weighted Euler sum formula})
\begin{equation}
\sum_{i=2}^{n-1} 2^i \zeta(i,n-i) = (n+1)\zeta (n),
\quad n\geq 2,
\mlabel{eq:oz}
\end{equation}
and applied it to study multiple zeta  star values.

In this paper we generalize the weighted Euler sum formula of Ohno and Zudilin to higher dimensions.
\begin{theorem}
{\bf (Weighted sum formula)} For positive integers $k\geq 2$ and $n\geq k+1$, we have
$$
\sum_{\substack{s_i\geq 1, s_1\geq 2\\
s_1+\cdots+s_{k}=n}}\hskip -10pt
\Big[2^{s_1-1}+(2^{s_1-1}-1)\Big(\big(\sum_{i=2}^{k-1}
2^{S_i-s_1-(i-1)}\big)+ 2^{S_{k-1}-s_1-(k-2)}\Big)\Big]
\zeta(s_1,\cdots, s_{k})=n\zeta(n),
$$
where $S_i=s_1+\cdots+s_i$ for $i=1,\cdots, k-1$.
\mlabel{thm:main}
\end{theorem}
See Theorem~\mref{thm:main1} for a more concise formulation. When $k=2$, this formula becomes
$$
\sum_{i=2}^{n-1} (2^i-1) \zeta(i,n-i) = n\zeta (n),
$$
which gives Eq.~(\mref{eq:oz}) after applying Eq.~(\mref{eq:euler}).

After introducing background notations and results, we prove our
weighted sum formula in Section~\mref{sec:main} by combining sum
formulas for the quasi-shuffle product, for the shuffle product, as
well as for \MZVs in Eq.~(\mref{eq:sum}). The proofs of the sum
formulas for the quasi-shuffle and shuffle product are given in
Section~\mref{sec:qsh} and Section~\mref{sec:sh} respectively. In
the two dimensional case~\mcite{OZ}, the sum formula for the shuffle
product was derived from Euler's decomposition formula. We have
generalized Euler's decomposition formula to multiple zeta
values~\mcite{GX2}. Instead of applying this generalization
directly, we obtain the sum formula for the shuffle product by
induction.

\medskip

\noindent {\bf Acknowledgements: }  Both authors thank the Max Planck Institute for Mathematics in Bonn for providing them with the environment to carried out this research. They also thank Don Zagier and Wadim Zudilin for helpful
discussions. The first author acknowledges the
support from NSF grant DMS-0505643.

\section{Double shuffle relations and the weighted sum theorem}
\mlabel{sec:main}
After introducing preliminary notations and results on MZVs in Section~\mref{ss:prel}, we prove the weighted sum formula (Theorem~\mref{thm:main}) in Section~\mref{ss:proof} by applying several shuffle and quasi-shuffle (stuffle) relations in Section~\mref{ss:proof}. The proofs of these relations will be given in the next two sections.
\subsection{Double shuffle relations of MZVs}
\mlabel{ss:prel} As is well-known, there are two ways to express a
\MZV:
\begin{eqnarray}
 \zeta(s_1,\cdots,s_k):&= &\sum_{n_1>\cdots>n_k\geq 1} \frac{1}{n_1^{s_1}\cdots n_k^{s_k}}
 \mlabel{eq:sumrep}\\
 &=&
 \int_0^1 \int_0^{t_1}\cdots \int_0^{t_{|\vec{s}|-1}} \frac{dt_1}{f_1(t_1)}
 \cdots \frac{dt_{|\vec{s}|}}{f_{|\vec{s}|}(t_{|\vec{s}|})}
 \mlabel{eq:int}
\end{eqnarray}
for integers $s_i\geq 1$ and $s_1>1$. Here
$$f_j(t)=\left\{\begin{array}{ll} 1-t_j, & j= s_1,s_1+s_2,\cdots, s_1+\cdots +s_k,\\
t_j, & \text{otherwise}. \end{array} \right .
$$
The product of two such sums is a $\ZZ$-linear combination of other
such sums and the product of two such integrals is a $\ZZ$-linear
combination of other such integrals. Thus the $\ZZ$-linear span of the \MZVs forms
an algebra
\begin{equation}
\mzvalg: = \ZZ \{ \zeta(s_1,\cdots,s_k)\ |\ s_i\geq 1, s_1\geq 2\}.
\mlabel{eq:mzvalg}
\end{equation}

The multiplication rules of the \MZVs reflected by these two
representations are encoded in two algebras, the quasi-shuffle (stuffle) algebra for the sum representation and the shuffle algebra for the integral representation~\mcite{Ho1,IKZ}.

For the sum representation, let $\calh\qsh$ be the quasi-shuffle
algebra whose underlying additive group is that of the
noncommutative polynomial algebra
$$ \ZZ\langle z_s\ |\ s\geq 1\rangle =\ZZ M (z_s\ |\ s\geq 1)$$
where $M(z_s\ |\ i\geq 1)$ is the free monoid generated by the set
$\{z_s\,|\, s\geq 1\}$. So an element of $M(z_s\ |\ s\geq 1)$ is
either the unit $1$, also called the empty word, or is of the form
$$z_{s_1,\cdots,s_k}:=z_{s_1}\cdots z_{s_k}, s_j\geq 1, 1\leq j\leq k, k\geq 1.$$
The product on $\calh\qsh$ is the
quasi-shuffle (also called harmonic shuffle or stuffle)
product~\mcite{3BL,Ho1,Ho2} defined recursively by
\begin{equation}
 (z_{r_1}u) \ast (z_{s_1} v) = z_{r_1}(u\ast (z_{s_1}v))
+ z_{s_1}((z_{r_1} u)\ast v) + z_{r_1+s_1}(u\ast v), \quad  u,v\in
M(z_s\ |\ s\geq 1) \mlabel{eq:qshde}
\end{equation}
with the convention that $1\ast u = u\ast 1= u$ for $u\in M(z_s\ | \
s\geq 1)$.
In the quasi-shuffle algebra $(\calh\qsh,\ast)$ there is a
subalgebra
\begin{equation}
\calh\qsh\lzero:=\ZZ \oplus \Big (\bigoplus_{s_1>1} \ZZ z_{s_1}
\cdots z_{s_k} \Big ) \subseteq \calh\qsh. \mlabel{eq:qssub}
\end{equation}
Then the multiplication rule of \MZVs according to their summation
representation follows from the statement that the linear map
\begin{equation}
\zeta\qsh: \calh\qsh\lzero \to \mzvalg, \quad z_{s_1,  \cdots, s_k} \mapsto \zeta(s_1,\cdots,s_k). \mlabel{eq:qshmzv}
\end{equation}
is an algebra homomorphism.

For the integral representation, let $\calh\shf$ be the shuffle
algebra whose underlying additive group is that of the
noncommutative polynomial algebra
$$\ZZ\langle x_0,x_1\rangle= \ZZ M (x_0,x_1)$$
where $M (x_0,x_1)$ is the free monoid generated by $x_0$ and $x_1$. The
product on $\calh\shf$ is the shuffle product defined recursively by
\begin{equation}
 (a  u) \ssha (b  v) = a  (u\ssha (b  v))+b  ((a  u)\ssha v), a,b\in \{x_0,x_1\}, u,v\in M (x_0,x_1)
 \mlabel{eq:shde}
\end{equation}
with the convention that $u \ssha 1=1\ssha u = u$ for $u\in M(x_0,x_1)$.
In the shuffle algebra $\calh\shf$, there are subalgebras
\begin{equation}
 \calh\shf\shzero:=\ZZ \oplus x_0\calh\shf x_1
\subseteq \calh\shf\lone:=\ZZ \oplus \calh\shf x_1 \subseteq
\calh\shf. \mlabel{eq:shsub}
\end{equation}
Then the multiplication rule of \MZVs according to their integral
representations follows from the statement that the linear map
\begin{equation}
 \zeta\shf: \calh\shf\shzero \to \mzvalg,
 \quad x_0^{  s_1-1}  x_1  \cdots   x_0^{s_k-1}  x_1 \mapsto \zeta(s_1,\cdots,s_k)
\mlabel{eq:shmzv}
\end{equation}
is an algebra homomorphism.

There is a natural bijection of abelian groups (but {\em not}
algebras)
\begin{equation}
\shqs: \calh\shf\lone \to \calh\qsh, \quad 1\leftrightarrow 1,\
x_0^{s_1-1}  x_1  \cdots   x_0^{s_k-1}  x_1
 \leftrightarrow z_{s_1, \cdots, s_k}.
\mlabel{eq:shqsh1}
\end{equation}
that restricts to a bijection of abelian groups
\begin{equation}
\shqs: \calh\shf\shzero \to \calh\qsh\lzero, \quad 1\leftrightarrow
1,\  x_0^{s_1-1}  x_1  \cdots   x_0^{s_k-1}  x_1
 \leftrightarrow z_{s_1, \cdots, s_k}.
\mlabel{eq:shqsh}
\end{equation}
Then the fact that \MZVs can be multiplied in two ways is reflected
by the commutative diagram
\begin{equation}
\xymatrix{ \calh\shf\shzero \ar^{\shqs}[rr] \ar_{\zeta\shf}[rd] &&
    \calh\qsh\lzero \ar^{\zeta\qsh}[ld]\\
& \mzvalg & } \mlabel{eq:shqshcom}
\end{equation}

Through $\shqs$, the shuffle product $\ssha$ on $\calh\shf\lone$ and
$\calh\shf\shzero$ transports to a product $\qssha$ on $\calh\qsh$
and $\calh\qsh\lzero$. That is, for $w_1,w_2\in \calh\qsh\lzero$,
define
\begin{equation}
w_1 \qssha w_2: = \shqs( \shqs^{-1}(w_1)\ssha \shqs^{-1}(w_2)).
\mlabel{eq:shtrans}
\end{equation}
Then the {\bf double shuffle relation} is simply the set
$$
\{ w_1\qssha w_2 - w_1 \ast w_2\ |\ w_1,w_2\in\calh\qsh\lzero\}$$
and the {\bf extended double shuffle relation}~\mcite{IKZ,Ra} is the set
\begin{equation}
\{ w_1\qssha w_2 - w_1 \ast w_2,\ z_1 \qssha w_2 - z_1 \ast w_2\ |\
w_1,w_2\in\calh\qsh\lzero\}.
\mlabel{eq:eds}
\end{equation}
\begin{theorem} {\bf (\mcite{Ho1,IKZ,Ra})}
Let $I_\edsalg$ be the ideal of $\calh\qsh\lzero$ generated by the extended double shuffle relation in Eq.~(\mref{eq:eds}). Then $I_\edsalg$ is in the kernel of $\zeta\qsh$.
\mlabel{thm:eds}
\end{theorem}
It is conjectured that $I_\edsalg$ is in fact the kernel of $\zeta\qsh$.

\subsection{Proof of the weighted sum theorem}
\mlabel{ss:proof}
Let $k$ be a positive integer $\geq 2$. For positive integers
$t_1,\cdots, t_{k-1}$, denote
\begin{equation}\begin{aligned} \calc(t_1,\cdots, t_{k-1}):
&=
(\sum_{j=1}^{k-1}2^{t_1+\cdots+t_j-j})+2^{t_1+\cdots+t_{k-1}-(k-1)}
\\
&=(\sum_{j=1}^{k-2}2^{t_1+\cdots+t_j-j})+
2^{t_1+\cdots+t_{k-1}-(k-1)+1}. \mlabel{eq:calc}
\end{aligned}\end{equation}
They satisfy the following simple relations that can be checked easily from the definition:
\begin{eqnarray}
\calc(t_1+1,t_2,\cdots, t_{k-1})= 2\calc(t_1, t_2, \cdots, t_{k-1}),
\mlabel{eq:coef1}
\\
\calc(1, t_2,\cdots, t_{k-1})= \calc(t_2,\cdots, t_{k-1})+1
\mlabel{eq:coef2}
\end{eqnarray} with the convention that $\calc(t_2,\cdots, t_{k-1})=
\calc(\emptyset)=1$ when $k=2$.

Using the notation $\calc(t_1,\cdots, t_{k-1})$ we can restate
Theorem \mref{thm:main} as follows.

\begin{theorem} {\bf (Second form of weighted sum formula)}
For positive integers $k$ and $n\geq k+1$, we
have
\begin{equation}
\sum_{\substack{ t_i\geq 1, t_1\geq 2 \\
t_1+\cdots +s_k=n}} [\calc(t_1,\cdots,
t_{k-1})-\calc(t_2,\cdots, t_{k-1})]\zeta(t_1, \cdots, t_{k})
=n\zeta(n),
\mlabel{eq:main1}
\end{equation}
with the convention that $\calc(t_2,\cdots, t_{k-1})=
\calc(\emptyset)=1$ when $k=2$. \mlabel{thm:main1}
\end{theorem}
We will prove Theorem~\mref{thm:main1} and hence the weighted sum
formula in Theorem~\mref{thm:main} by using the following auxiliary
Theorem~\mref{thm:qsh1} and Theorem~\mref{thm:sha1}. They are a type
of sum formula on the products $\ast$ and $\qssha$ respectively and
are interesting on their own right. Their proofs will be given in
Section~\mref{sec:qsh} and Section~\mref{sec:sh} respectively.

We first display the sum formulas on the product $*$.

\begin{theorem}
For positive integers $k\geq 2$ and $n\geq k$, we have
\begin{equation}
\sum_{\substack{ r, s_i\geq 1 \\ r+s_1+\cdots +s_{k-1}=n
}} \hskip -10pt z_{r}
* z_{s_1,\cdots, s_{k-1}} = k \hskip -10pt
\sum_{\substack{t_i\geq 1 \\ t_1+\cdots+t_{k}=n
}} \hskip -10pt  z_{t_1,\cdots, t_{k}} \ +
\ (n-k+1)\hskip -20pt \sum_{\substack{u_i\geq 1 \\
u_1+\cdots+u_{k-1}=n }} \hskip -20pt z_{u_1,\cdots,
u_{k-1}}. \mlabel{eq:qsh}
\end{equation}
\mlabel{thm:qsh}
\end{theorem}
This theorem will be applied to \MZVs through the following
\begin{theorem} For positive integers $k\geq 2$ and  $n\geq k+1,$ we have
\begin{equation} \begin{aligned}
& \sum_{\substack{ r, s_i\geq 1, s_1\geq 2 \\
r+s_1+\cdots +s_{k-1}=n
}}  z_{r}
* z_{s_1,\cdots, s_{k-1}}
\\
= & \sum_{\substack{ t_i\geq 1, t_1=1,t_2\geq 2 \\
t_1+\cdots+t_{k}=n}} z_{t_1,t_2,\cdots, t_{k}} +(k-1)
\hskip -5pt
\sum_{\substack{t_i\geq 1, t_1\geq 2, t_2=1 \\
t_1+\cdots+t_{k}=n }} \hskip -10pt z_{t_1,\cdots,t_{k}}
 \\
 & \quad \quad +k
\hskip -10pt
\sum_{\substack{t_i\geq 1, t_1\geq 2, t_2\geq 2 \\
t_1+\cdots+t_{k}=n
}} \hskip -10pt  z_{t_1,\cdots, t_{k}} \ +
\ (n-k)\hskip -20pt \sum_{\substack{u_i\geq 1, u_1\geq 2 \\
u_1+\cdots+u_{k-1}=n }} \hskip -20pt z_{u_1,\cdots,
u_{k-1}}.
\end{aligned} \mlabel{eq:qsh1}
\end{equation}
\mlabel{thm:qsh1}
\end{theorem}
The proofs of these two theorems will be given in Section~\mref{sec:qsh}.
Similarly for the product $\qssha$, we have the following sum formulas.
\begin{theorem} For positive integers $k\geq 2$ and $n\geq k$, we have
\begin{equation}
\sum_{r, s_i\geq 1, r+s_1+\cdots +s_{k-1}=n} z_{r}\qssha
z_{s_1,\cdots, s_{k-1}} = \sum_{t_i\geq 1, t_1+\cdots+t_{k}=n}
\calc(t_1,\cdots, t_{k-1}) z_{t_1,\cdots, t_{k}}. \mlabel{eq:sha}
\end{equation} \mlabel{thm:sha}
\end{theorem}

\begin{theorem} For positive integers $k\geq 2$ and $n\geq k+1$, we have
{\allowdisplaybreaks
\begin{equation}
\begin{aligned} & \sum_{\substack{r, s_i\geq 1, s_1\geq 2\\
r+s_1+\cdots +s_{k-1}=n}}\hskip -10pt  z_{r}\qssha
z_{s_1,\cdots, s_{k-1}} \\
=& \hskip -10pt \sum_{
\substack{t_i\geq 1\\
t_1+\cdots+t_{k}=n }}\hskip -10pt [ \calc(t_1,\cdots,
t_{k-1})- \calc(t_2,\cdots, t_{k-1})] z_{t_1,\cdots,
t_{k}}  - \hskip -10pt \sum_{ \substack{ t_i\geq 1, t_2=1
\\
t_1+\cdots+t_{k}=n }}\hskip -10pt  z_{t_1,\cdots, t_{k}}
\mlabel{eq:sha1}
\end{aligned}\end{equation}}
with the convention that $\calc(t_2,\cdots, t_{k-1})=
\calc(\emptyset)=1$ when $k=2$. \mlabel{thm:sha1}
\end{theorem}
The proofs of these two theorems will be given in Section~\mref{sec:sh}.
Now we derive Theorem \mref{thm:main1} and hence Theorem~\mref{thm:main} from Theorem \mref{thm:qsh1} and
Theorem \mref{thm:sha1}. \vskip
6pt

\noindent{\it Proof of Theorem~\mref{thm:main1}.} {} Regrouping the
sums in Eq.~(\mref{eq:sha1}) of Theorem \mref{thm:qsh1} and applying
the summation relation $$\sum\limits_{\substack{t_1=1, t_2\geq 2\\
t_1+\cdots+t_k=n}}-\sum\limits_{\substack{t\geq 2, t_2=1\\
t_1+\cdots+t_k=n}}=\sum\limits_{\substack{t_1=1\\
t_1+\cdots +t_k=n}}-\sum\limits_{\substack{t_2=1 \\ t_1+\cdots+t_k=n
}}$$ we obtain
\begin{eqnarray*}
&& \sum_{\substack{ r, s_i\geq 1, s_1\geq 2 \\
r+s_1+\cdots +s_{k-1}=n
}}z_r*z_{s_1,\cdots, s_{k-1}} \\
&=& \hspace{-.3cm} \sum_{\substack{ t_i\geq
1, t_1=1 \\ t_1+\cdots+t_k=n }} z_{t_1,\cdots, t_k}+
k \hspace{-.3cm}\sum_{\substack{t_i\geq 1, t_1\geq 2 \\
t_1+\cdots+t_k=n}} z_{t_1,\cdots, t_k}
 -\hspace{-.3cm}\sum_{\substack{t_i\geq 1, t_2=1\\ t_1+\cdots+t_k=n
}}z_{t_1,\cdots, t_k} +(n-k)\hspace{-.3cm}
\sum_{\substack{ u_i\geq 1, u_1\geq 2 \\
u_1+\cdots+u_{k-1}=n }}z_{u_1,\cdots, u_{k-1}}.
\end{eqnarray*}
From this equation and Theorem \mref{thm:sha1} we obtain
{\allowdisplaybreaks
\begin{eqnarray*}
&&\sum_{\substack{ r, s_i\geq 1, s_1\geq 2 \\
r+s_1+\cdots +s_{k-1}=n}} (z_r\qssha z_{s_1,\cdots,
s_{k-1}}- z_r*z_{s_1,\cdots, s_{k-1}})
\\
&=&\sum_{\substack{t_i\geq 1, t_1=1\\
t_1+\cdots+t_k=n}}[\calc(t_1,\cdots,
t_{k-1})-\calc(t_2,\cdots, t_{k-1})-1]z_{t_1,\cdots, t_k}
\\ && \quad \quad +
\sum_{\substack{ t_i\geq 1, t_1\geq
2\\
t_1+\cdots+t_k=n }} [\calc(t_1,\cdots,
t_{k-1})-\calc(t_2,\cdots, t_{k-1})-k]z_{t_1,\cdots, t_k}
\\
&& \quad \quad -(n-k)\sum_{\substack{u_i\geq 1, u_1\geq 2 \\
u_1+\cdots+u_{k-1}=n }}z_{u_1,\cdots, u_{k-1}}
\\
&= &\sum_{\substack{ t_i\geq 1, t_1\geq
2\\
t_1+\cdots+t_k=n }} [\calc(t_1,\cdots, t_{k-1})-\calc(t_2,\cdots,
t_{k-1})-k]z_{t_1,\cdots, t_k}
\\
&& \quad \quad -(n-k)\sum_{\substack{u_i\geq 1, u_1\geq 2 \\
u_1+\cdots+u_{k-1}=n }}z_{u_1,\cdots, u_{k-1}} .
\end{eqnarray*}
} by Eq.~(\mref{eq:coef2}). Since $\zeta^*(z_r\qssha z_{s_1,\cdots,
s_{k-1}}- z_r*z_{s_1,\cdots, s_{k-1}})=0$, by
Theorem~\mref{thm:eds}, this gives
$$\begin{aligned} & \sum_{\substack{ t_i\geq 1, t_1\geq 2 \\
t_1+\cdots +t_k=n}} [\calc(t_1,\cdots,
t_{k-1})-\calc(t_2,\cdots, t_{k-1})-k]\zeta(t_1, \cdots, t_{k})
\\
&= (n-k)\sum_{\substack{u_i\geq 1, u_1\geq 2 \\
u_1+\cdots+u_{k-1}=n }}\zeta(u_1,\cdots, u_{k-1}).
\end{aligned}$$ The sum formula in Eq.~(\mref{eq:sum}) shows that
$$\sum_{\substack{t_i\geq 1, t_1\geq 2 \\
t_1+\cdots +t_k=n}} \zeta(t_1, \cdots, t_{k})
=\sum_{\substack{u_i\geq 1, u_1\geq 2 \\
u_1+\cdots+u_{k-1}=n }}\zeta(u_1,\cdots,
u_{k-1})=\zeta(n).$$ Therefore
$$\sum_{\substack{ t_i\geq 1, t_1\geq 2 \\
t_1+\cdots +t_k=n}} [\calc(t_1,\cdots,
t_{k-1})-\calc(t_2,\cdots, t_{k-1})]\zeta(t_1, \cdots, t_{k})
=n\zeta(n),
$$ as desired.
\qed

\section{Proofs of Theorem \mref{thm:qsh} and Theorem
\mref{thm:qsh1}}
\mlabel{sec:qsh}
In this section we prove the two sum relations of the quasi-shuffle product.

\subsection{Proof of Theorem \mref{thm:qsh}}
We first consider the case when $k=2$. By the basic relation
$z_r\ast z_{s_1}=z_{r,s_1}+z_{s_1,r}+z_{r+s_1}$ from
Eq.~(\mref{eq:qshde}), we have
\begin{eqnarray*}
\sum_{r,s_1\geq 1, r+s_1=n} (z_r\ast z_{s_1})&=&\sum_{r,s_1\geq 1, r+s_1=n}z_{r,s_1}+\sum_{r,s_1\geq 1, r+s_1=n}z_{s_1,r}+\sum_{r,s_1\geq 1, r+s_1=n}z_{r+s_1}\\
&=& 2\sum_{t_i\geq 1, t_1+t_2=n} z_{t_1,t_2} + (n-1)z_{n},
\end{eqnarray*}
as needed.

In the general case we prove Eq.~(\mref{eq:qsh}) by induction on
$n$. If $n=2$, then $k=2$ and we are done.
For a given integer $m\geq 3$, assume that Eq.~(\mref{eq:qsh})
holds when $n<m$ and consider Eq.~(\mref{eq:qsh}) when $n=m$ and $k\geq 2$. Since we have proved the case when $k=2$, we may assume that $k\geq 3$.
By Eq.~(\mref{eq:qshde}),
we have
\begin{equation}
\begin{aligned}
&\sum_{\substack{ r, s_i\geq 1\\ r+s_1+\cdots +s_{k-1}=n}} (z_r* z_{s_1,\cdots, s_{k-1}})
\\
& =
\hskip -20pt \sum_{\substack{r, s_i\geq 1\\
r+s_1+\cdots+s_{k-1}=n}}\hskip -20pt  z_{r,s_1,\cdots,
s_{k-1}} +
\hskip -20pt  \sum_{\substack{ r, s_i\geq 1\\
 r+s_1+\cdots+s_{k-1}=n}} \hskip -20pt
z_{s_1}(z_r* z_{s_2,\cdots, s_{k-1}})
+\hskip -20pt
 \sum_{\substack{ r, s_i\geq 1\\
r+s_1+\cdots+s_{k-1}=n}} \hskip -20pt z_{r+s_1,s_2,\cdots,
s_{k-1}}.
\end{aligned}
\label{eq:qshrec}
\end{equation}
In the second sum on the right hand side of Eq.~(\mref{eq:qshrec}),
for any fixed $s_1\geq 1$, by the induction hypothesis we have
$$ \begin{aligned} & \sum_{r,s_i\geq 1, r+s_2+\cdots+s_{k-1}=n-s_1}z_r* z_{s_2,\cdots, s_{k-1}}
\\
& = (k-1)\sum_{t_i\geq 1, t_1+\cdots+t_{k-1}=n-s_1}z_{t_1,\cdots,
t_{k-1}} + (n-s_1-k+2)\sum_{u_i\geq
1,u_1+\cdots+u_{k-2}=n-s_1}z_{u_1,\cdots, u_{k-2}}.
\end{aligned} $$ So the second sum becomes
{\allowdisplaybreaks
\begin{eqnarray*}
& & \sum_{\substack{r,s_i\geq 1\\r+s_1+s_2+\cdots+s_{k-1}=n}}z_{s_1}(z_r* z_{s_2,\cdots,
s_{k-1}})
\\
& = &(k-1) \hspace{-.6cm}
\sum_{\substack{ s_1\geq 1, t_i\geq 1 \\
s_1+t_1+\cdots+t_{k-1}=n }}  z_{s_1,t_1,\cdots, t_{k-1}} +
(n-s_1-k+2) \hspace{-.6cm}\sum_{\substack{ s_1\geq 1, u_i\geq
1 \\ s_1+u_1+\cdots+u_{k-2}=n }}  z_{s_1,u_1,\cdots, u_{k-2}} \\
&= &(k-1)\sum_{t_i\geq 1, t_1+\cdots+t_k=n}  z_{t_1,\cdots, t_{k}} +
(n-u_1-k+2)\sum_{u_i\geq 1,u_1+\cdots+u_{k-1}=n}z_{u_1,\cdots,
u_{k-1}}.
\end{eqnarray*}
}

\noindent For the third sum on the right hand side of
Eq.~(\mref{eq:qshrec}), we have
$$\sum_{r, s_i\geq 1,r+s_1+\cdots +s_{k-1}=n}
z_{r+s_1,s_2,\cdots, s_{k-1}} = \sum_{u_i\geq 1,
u_1+\cdots+u_{k-1}=n}(u_1-1) z_{u_1,\cdots, u_{k-1}} .$$
Therefore Eq.~(\mref{eq:qshrec}) becomes
{\allowdisplaybreaks
\begin{eqnarray*}
& & \sum_{\substack{r, s_i\geq 1\\ r+s_1+\cdots+s_{k-1}=n}} z_r* z_{s_1,\cdots, s_{k-1}}
\\
&=&  \sum_{t_i\geq 1, t_1+\cdots+t_{k}=n}z_{t_1,\cdots, t_{k}} + \sum_{t_i\geq 1, t_1+\cdots+t_k=n} (k-1)
   z_{t_1,\cdots, t_{k}}\\
   &&  + \sum_{u_i\geq
    1,u_1+\cdots+u_{k-1}=n}(n-u_1-k+2)z_{u_1,\cdots, u_{k-1}}
    +\sum_{u_i\geq 1, u_1+\cdots+u_{k-1}=n}(u_1-1)
     z_{u_1,\cdots, u_{k-1}}
 \\
& = &\sum_{t_i\geq 1, t_1+\cdots+t_k=n} k
   z_{t_1,\cdots, t_{k}} +\sum_{u_i\geq
    1,u_1+\cdots+u_{k-1}=n}(n-k+1)z_{u_1,\cdots, u_{k-1}}.
\end{eqnarray*}
}

\noindent This means that Eq.~(\mref{eq:qsh}) holds when $n=m$,
completing the inductive proof of Theorem~\mref{thm:qsh}.

\subsection{Proof of Theorem \mref{thm:qsh1}}
We now prove Theorem~\mref{thm:qsh1} by applying
Theorem~\mref{thm:qsh}. When $k=2$ we can verify
Eq.~(\mref{eq:qsh1}) directly by Eq.~(\mref{eq:qshde}) as in the case of Theorem~\mref{thm:qsh}. For $k\geq
3$, applying Eq.~(\mref{eq:qshde}) and Eq.~(\mref{eq:qsh}), we have
{\allowdisplaybreaks
\begin{eqnarray*}
&& \sum_{\substack{ r,s_i\geq 1 \\
r+s_2+\cdots+s_{k-1}=n-1
}} z_r* z_{1,s_2,\cdots, s_{k-1}}
\\
& =&\hskip -10pt
\hskip -10pt \sum_{\substack{ r,s_i\geq 1 \\
r+s_2+\cdots+s_{k-1}=n-1 }} \hskip -10pt z_{r, 1, s_2, \cdots,
s_{k-1}}
 +
 \sum_{\substack{ r,s_i\geq 1 \\
r+s_2+\cdots+s_{k-1}=n-1 }} \hskip -10pt z_1(z_r * z_{s_2,\cdots,
s_{k-1}}) + \sum_{\substack{ r,s_i\geq 1 \\
r+s_2+\cdots+s_{k-1}=n-1 }} \hskip -10pt  z_{r+1,s_2,\cdots,
s_{k-1}}
\\
& =& \sum_{\substack{t_i\geq 1, t_2=1\\ t_1+\cdots +t_{k}=n}}
z_{t_1,\cdots, t_{k}}
  + \bigg((k-1) \sum_{t_2+\cdots+t_{k}=n-1}z_{1,t_2,\cdots, t_{k}
}
\\
&& +[(n-1)-(k-1)+1]\sum_{u_2+\cdots+u_{k-1}=n-1}z_{1, u_2, \cdots,
u_{k-1}}\bigg) + \sum_{\substack{u_i\geq 1, u_1\geq 2\\
u_1+\cdots+u_{k-1}=n}}z_{u_1,\cdots, u_{k-1}},
\end{eqnarray*}
}

\noindent where $k$ and $n$ in Eq.~(\mref{eq:qsh}) are replaced by
$k-1$ and $n-1$.  By regrouping, this expression can be further
simplified to {\allowdisplaybreaks
\begin{eqnarray*}
&& (k-1)\sum_{\substack{t_i\geq 1,t_1=1\\
t_1+\cdots+t_{k}=n}}z_{t_1, t_2,\cdots,t_{k}} +
\sum_{\substack{t_i\geq 1, t_2=1\\ t_1+\cdots+t_{k}=n}}z_{t_1,
t_2,\cdots, t_{k}}
\\ & & \hskip 20pt
+ (n-k+1)\sum_{\substack{u_i\geq 1,u_1=1\\
u_1+\cdots+u_{k-1}=n}}z_{u_1, u_2,\cdots, u_{k-1}} +
\sum_{\substack{u_i\geq 1,u_1\geq 2\\ u_1+\cdots+u_{k-1}=n}}z_{u_1,
u_2,\cdots, u_{k-1}}
\\
& = & (k-1)\sum_{\substack{t_i\geq 1, t_1=1, t_2\geq 2\\
t_1+\cdots+t_{k}=n}}z_{t_1, t_2,\cdots,t_{k}} +
\sum_{\substack{t_i\geq 1, t_1\geq 2, t_2=1\\
t_1+\cdots+t_{k}=n}}z_{t_1, t_2,\cdots, t_{k}}
\\
&&  + k \sum_{\substack{t_i\geq 1, t_1=t_2=1\\
t_1+\cdots+t_{k}=n}}z_{t_1, t_2,\cdots, t_{k}}
 + (n-k+1)\sum_{\substack{u_i\geq 1,u_1=1\\ u_1+\cdots+u_{k-1}=n}}z_{u_1,
u_2,\cdots, u_{k-1}}
\\
&& + \sum_{\substack{u_i\geq 1,u_1\geq 2\\
u_1+\cdots+u_{k-1}=n}}z_{u_1, u_2,\cdots, u_{k-1}}.
\end{eqnarray*}
}

\noindent Now Eq.~(\mref{eq:qsh1}) follows from this and Eq.
(\mref{eq:qsh}).

\section{Proofs of Theorem \mref{thm:sha} and Theorem
\mref{thm:sha1}} \mlabel{sec:sh}
This sections gives the proofs of the two sum relations of the shuffle product.
\subsection{A preparational lemma} In this subsection we provide a
lemma that is needed in the proofs of Theorem \mref{thm:sha} and
Theorem \mref{thm:sha1}.

Let $\calh\qsh^+$ be the subring of $\calh\qsh$ generated by
$z_{\vec{s}}$ with $\vec{s}\in \ZZ_{\geq 1}^k, k\geq 1$. Then
$$ \calh\qsh= \ZZ\oplus \calh\qsh^+. $$
Define two operators
\begin{equation}
\begin{aligned}
P:& \calh\qsh^+\rightarrow \calh\qsh,
P(z_{s_1,s_2,\cdots,s_k})= z_{s_1+1,s_2,\cdots, s_k}, \\
Q:& \calh\qsh\rightarrow \calh\qsh, Q(w)=\left \{\begin{array}{ll}
z_1w, & w\neq 1, \\ z_1, & w=1.
\end{array} \right .
\end{aligned}
\mlabel{eq:pq}
\end{equation}
These operators are simply the transports of the operators $$
\begin{aligned}
I_0: & \calh\shf\lone^+(\ssg{\sg})\to \calh\shf\lone(\ssg{\sg}), \quad I_0(u)= x_0 u,\\
I_1: & \calh\shf\lone(\ssg{\sg}) \to \calh\shf\lone(\ssg{\sg}),
\quad I_1(u)=\left \{\begin{array}{ll} x_1 u, & u \neq 1, \\
x_1, & u=1.
\end{array}\right .
\end{aligned}
$$
Thus the recursive definition of $\ssha$ in Eq.~(\mref{eq:shde}) gives the following Rota-Baxter relation~\mcite{Ba,GK1,GX1,Ro} from~\mcite{GX2}.
\begin{prop} $($\mcite{GX2}$)$ The multiplication $\qssha$ on
$\calh\qsh$ defined in Eq.~(\mref{eq:shtrans}) is the unique one
such that
\begin{equation}
\begin{aligned}
P(w_1)\qssha P(w_2) &= P\big(w_1 \qssha P(w_2)\big)+
P\big(P(w_1)\qssha w_2\big), \qquad w_1,w_2\in\calh\qsh^+,
\\
Q(w_1)\qssha Q(w_2) &= Q\big(w_1 \qssha Q(w_2)\big)+ Q\big(
Q(w_1)\qssha w_2\big), \qquad w_1,w_2\in\calh\qsh,
\\
P(w_1)\qssha Q(w_2) &= Q\big(P(w_1)\qssha w_2\big)+ P\big(w_1\qssha
Q(w_2)\big), \qquad w_1\in\calh\qsh^+, w_2 \in\calh\qsh,
\\
Q(w_1)\qssha P(w_2) &= Q\big(w_1\qssha P(w_2)\big)+
P\big(Q(w_1)\qssha w_2\big), \qquad w_1\in\calh\qsh, w_2
\in\calh\qsh^+.
\end{aligned}
\mlabel{eq:iqform}
\end{equation}
with the initial condition that $1 \qssha w=w\qssha 1 = w$ for $w\in
\calh\qsh$. \mlabel{pp:rec}
\end{prop}

\begin{lemma} For positive integers $k\geq 2$ and $n\geq k$,
we have
\begin{equation}
\begin{aligned}& z_1\qssha  \sum_{\substack{s_i\geq 1, s_1\geq 2\\ s_1+\cdots+s_{k-1}=n-1}}z_{s_1, \cdots, s_{k-1}}
\\ &= \sum_{\substack{t_i\geq 1,t_1=1, t_2\geq 2\\
t_1+\cdots+t_{k}=n}}z_{t_1, t_2,\cdots, t_{k}}+ k \hskip -2pt
\sum_{\substack{ t_i\geq 1, t_1\geq 2, t_k=1 \\
t_1+\cdots+t_{k}=n}} \hskip -5pt z_{t_1,\cdots, t_{k}} + (k-1)
\hskip -5pt \sum_{\substack{ t_i\geq 1, t_1, t_k\geq 2 \\
t_1+\cdots+t_{k}=n}}\hskip -5pt z_{t_1,\cdots, t_{k}}.
\end{aligned}
\mlabel{eq:z1b}
\end{equation}
 and
\begin{equation}
z_1\qssha \hskip -15pt \sum_{\substack{ s_i\geq 1\\
s_1+\cdots+s_{k-1}=n-1}} \hskip -15pt  z_{s_1, \cdots, s_{k-1}}= k
\hskip -10pt \sum_{\substack{ t_i\geq 1, t_k=1\\
t_1+\cdots+t_{k}=n}}z_{t_1,\cdots, t_{k}} + (k-1)\hskip -10pt
\sum_{\substack{ t_i\geq 1, t_k\geq 2 \\
t_1+\cdots+t_{k}=n }} z_{t_1,\cdots, t_{k}}. \mlabel{eq:z1a}
\end{equation}
\mlabel{lem:z1sha}
\end{lemma}
Both sides of Eq.~(\mref{eq:z1b}) are zero when $n=k$.

\begin{proof} We prove Eq.~(\mref{eq:z1b}) and (\mref{eq:z1a}) by
induction on $n$. If $n=2, 3$, Eq.~(\mref{eq:z1b}) and
(\mref{eq:z1a}) can be verified directly. Let $m\geq 4$ be an
integer. Assume that Eq.~(\mref{eq:z1b}) and (\mref{eq:z1a}) hold
when $n< m$. Now assume that $n=m$.

By Eq.~(\mref{eq:pq}), Eq.~(\mref{eq:iqform}) and the induction hypothesis, the left hand side of Eq.~(\mref{eq:z1b}) becomes
{\allowdisplaybreaks
\begin{eqnarray*}
\lefteqn{z_1\qssha \hskip -10pt \sum_{\substack{s_i\geq 1, s_1\geq 2\\
s_1+\cdots+s_{k-1}=n-1}} \hskip -15pt z_{s_1,\cdots,
s_{k-1}}
 = Q(1)\qssha P(\sum_{s_i\geq 1, s_1+\cdots+s_{k-1}=n-2} z_{s_1, \cdots, s_{k-1}})}\\
& = &Q(\hskip -11pt \sum_{\substack{s_i\geq 1,s_1\geq 2\\
s_1+\cdots+s_{k-1}=n-1 } }\hskip -20pt  z_{s_1, \cdots,
s_{k-1}}) + P(z_1\qssha
\hskip -15pt \sum_{\substack{s_i\geq 1\\
 s_1+\cdots+s_{k-1}=n-2}} \hskip -15pt  z_{s_1,
\cdots, s_{k-1}})\\
&= &\sum_{\substack{t_i\geq 1, t_2\geq 2\\
t_2+\cdots+t_{k+1}=n-1}}z_{1, t_2,\cdots, t_{k+1}} +
P\Big(k\sum_{\substack{t_i\geq 1, t_k=1\\ t_1+\cdots +
t_{k}=n-1}}z_{t_1,\cdots, t_{k}}+(k-1)\sum_{\substack{t_i\geq 1, t_k\geq 2 \\
t_1+\cdots + t_{k}=n-1}}z_{t_1,\cdots, t_{k}}\Big)
\\
& =& \sum_{\substack{t_i\geq 1, t_2\geq 2\\
t_2+\cdots+t_{k+1}=n-1}}z_{1, t_2,\cdots, t_{k+1}} +
k\sum_{\substack{t_i\geq 1, t_1\geq 2,
t_k=1\\t_1+\cdots+t_{k}=n}}z_{t_1,\cdots, t_{k}}+
(k-1)\sum_{\substack{t_i\geq 1, t_1,t_k\geq
2\\t_1+\cdots+t_{k}=n}}z_{t_1,\cdots, t_{k}},
\end{eqnarray*}
}

\noindent which shows that Eq.~(\mref{eq:z1b}) holds when $n=m$.

By a similar argument, we have
{\allowdisplaybreaks
\begin{eqnarray*}
\lefteqn{ z_1\qssha \sum_{s_2+\cdots+s_{k-1}=n-2} z_{1, s_2,\cdots, s_{k-1}}}
   \\
 &=& Q(\sum_{s_2+\cdots+s_{k-1}=n-2} z_{1, s_2, \cdots, s_{k-1}})+
      Q(z_1\qssha \sum_{s_2+\cdots+s_{k-1}=n-2} z_{s_2,\cdots, s_{k-1}})
   \\
   & =& \sum_{s_i\geq 1, s_2+\cdots+s_{k-1}=n-2}z_{1,1,s_2,\cdots, s_{k-1}}
   \\
   &&+
   Q\Big((k-1) \sum_{\substack{t_i\geq 1,t_{k-1}=1 \\
   t_1+\cdots+t_{k-1}=n-1}}
   z_{t_1,\cdots,t_{k-1}}+ (k-2) \sum_{\substack{t_i\geq 1,t_{k-1} \geq 2 \\
   t_1+\cdots+t_{k-1}=n-1}}
   z_{t_1,\cdots,t_{k-1}}\Big)
   \\
   & = &\sum_{\substack{t_i\geq 1,t_1=t_2=1\\
    t_1+\cdots+t_{k}=n}}z_{t_1,t_2,\cdots, t_{k}} +
   (k-1) \sum_{\substack{t_i\geq 1,t_1=t_k=1 \\ t_1+\cdots+t_{k}=n}} z_{t_1, t_2,\cdots,t_{k}}
   +(k-2) \sum_{\substack{t_i\geq 1,t_1=1, t_k\geq 2 \\ t_1+\cdots+t_{k}=n}} z_{t_1, t_2,\cdots,t_{k}}
\end{eqnarray*}
}

\noindent Adding this with Eq.~(\mref{eq:z1b}) we obtain
{\allowdisplaybreaks
\begin{eqnarray*}
&& z_1\qssha \hskip -15pt \sum_{\substack{ s_i\geq 1\\
s_1+\cdots+s_{k-1}=n-1}} \hskip -15pt  z_{s_1, \cdots, s_{k-1}}
\\
&=&\sum_{\substack{t_i\geq 1,t_1=1, t_2\geq 2\\
t_1+\cdots+t_{k}=n}}z_{t_1, t_2,\cdots, t_{k}}+ k \hskip -2pt
\sum_{\substack{ t_i\geq 1, t_1\geq 2, t_k=1 \\
t_1+\cdots+t_{k}=n}} \hskip -5pt z_{t_1,\cdots, t_{k}} + (k-1)
\hskip -5pt \sum_{\substack{ t_i\geq 1, t_1, t_k\geq 2 \\
t_1+\cdots+t_{k}=n}}\hskip -5pt z_{t_1,\cdots, t_{k}}
\\
&& + \sum_{\substack{t_i\geq 1,t_1=t_2=1\\
    t_1+\cdots+t_{k}=n}}z_{t_1,t_2,\cdots, t_{k}} +
   (k-1) \sum_{\substack{t_i\geq 1,t_1=t_k=1 \\ t_1+\cdots+t_{k}=n}} z_{t_1, t_2,\cdots,t_{k}}
   +(k-2) \sum_{\substack{t_i\geq 1,t_1=1, t_k\geq 2 \\ t_1+\cdots+t_{k}=n}} z_{t_1, t_2,\cdots,t_{k}}
\\
&=& k \hskip -2pt
\sum_{\substack{ t_i\geq 1, t_1\geq 2, t_k=1 \\
t_1+\cdots+t_{k}=n}} \hskip -5pt z_{t_1,\cdots, t_{k}} + (k-1)
\hskip -5pt \sum_{\substack{ t_i\geq 1, t_1, t_k\geq 2 \\
t_1+\cdots+t_{k}=n}}\hskip -5pt z_{t_1,\cdots, t_{k}}
\\
&& + \sum_{\substack{t_i\geq 1,t_1=1\\
    t_1+\cdots+t_{k}=n}}z_{t_1,t_2,\cdots, t_{k}} +
   (k-1) \sum_{\substack{t_i\geq 1,t_1=t_k=1 \\ t_1+\cdots+t_{k}=n}} z_{t_1, t_2,\cdots,t_{k}}
   +(k-2) \sum_{\substack{t_i\geq 1,t_1=1, t_k\geq 2 \\ t_1+\cdots+t_{k}=n}} z_{t_1, t_2,\cdots,t_{k}}
\\&=& k \hskip -2pt
\sum_{\substack{ t_i\geq 1, t_1\geq 2, t_k=1 \\
t_1+\cdots+t_{k}=n}} \hskip -5pt z_{t_1,\cdots, t_{k}} + (k-1)
\hskip -5pt \sum_{\substack{ t_i\geq 1, t_1, t_k\geq 2 \\
t_1+\cdots+t_{k}=n}}\hskip -5pt z_{t_1,\cdots, t_{k}}
\\
&& +
   k\sum_{\substack{t_i\geq 1,t_1=t_k=1 \\ t_1+\cdots+t_{k}=n}} z_{t_1, t_2,\cdots,t_{k}}
   +(k-1) \sum_{\substack{t_i\geq 1,t_1=1, t_k\geq 2 \\ t_1+\cdots+t_{k}=n}} z_{t_1,
   t_2,\cdots,t_{k}}
\\&=& k \hskip -2pt
\sum_{\substack{ t_i\geq 1,  t_k=1 \\
t_1+\cdots+t_{k}=n}} \hskip -5pt z_{t_1,\cdots, t_{k}} + (k-1)
\hskip -5pt \sum_{\substack{ t_i\geq 1, t_k\geq 2 \\
t_1+\cdots+t_{k}=n}}\hskip -5pt z_{t_1,\cdots, t_{k}},
\end{eqnarray*}
}

\noindent which shows that Eq.~(\mref{eq:z1a}) holds when $n=m$.
This completes the induction.
\end{proof}

\subsection{Proof of Theorem \mref{thm:sha}}
For the proof of Theorem \mref{thm:sha}, we first consider
the case of $k=2$. In this case we have Euler's decomposition formula~\mcite{GKZ,GX2}
$$z_r\qssha z_{s_1}=
\sum_{t_1\geq 1, t_2\geq 1,
t_1+t_2=n}\Big(\binc{t_1-1}{r-1}+\binc{t_1-1}{s_1-1}\Big)
z_{t_1,t_2}.
$$
So
$$ \sum_{r\geq 1, s_1\geq 1, r+s_1=n}z_r\qssha z_{s_1}
 = \sum_{t_1\geq 1, t_2\geq 1,
t_1+t_2=n}\Big(
 \sum_{r= 1}^{t_1}
\binc{t_1-1}{r-1}+
 \sum_{s_1=1}^{t_1}
\binc{t_1-1}{s_1-1}\Big)
z_{t_1,t_2}
=   2^{t_1}z_{t_1,t_2},
$$
 as required.

For the general case we prove Eq.~(\mref{eq:sha}) by induction on
$n$. If $n=2$, then $k=2$ and so we are done by the above argument. Let $m$ be a
positive integer $\geq 3$ and assume that Eq.~(\mref{eq:sha}) holds
for $n<m$. Now assume that $n=m$. Since we have dealt with the case
of $k=2$, we may assume that $k\geq 3$.

We  decompose the left hands side of Eq.~(\mref{eq:sha}) into three
disjoint parts when $r=1, s_1 \geq 1$, when $r\geq 2, s_1\geq 2$ and
when $r\geq 2, s_1=1$: {\allowdisplaybreaks
\begin{equation}
\begin{aligned}
    \sum_{\substack{r, s_i\geq 1\\ r+s_1+\cdots +s_{k-1}=n}} z_{r}\qssha z_{s_1,\cdots, s_{k-1}}
= &
    \sum_{ \substack{s_i\geq 1\\ s_1+\cdots +s_{k-1}=n-1}} z_{1}\qssha z_{s_1,\cdots,
    s_{k-1}}
\\ & +
    \sum_{\substack{s_i\geq 1, r,s_1\geq 2\\  r+s_1+\cdots+s_{k-1}=n}}z_r \qssha
    z_{s_1,\cdots,s_{k-1}}
\\ & +
    \sum_{\substack{r\geq 2, s_i\geq 1\\
    r+s_2+\cdots+s_{k-1}=n-1}}z_r\qssha z_{1,s_2,\cdots, s_{k-1}}.
\end{aligned}
\mlabel{eq:shsum}
\end{equation}}
We denote the three sums by $\sumx, \sumy$ and $\sumz$ respectively with the given order.

For the sum $\sumy$, by Eq.~(\mref{eq:pq}) and (\mref{eq:iqform}),
we have {\allowdisplaybreaks
\begin{eqnarray*}
&& \sum_{\substack{s_i\geq 1, r,s_1\geq 2\\
r+s_1+\cdots+s_{k-1}=n}}z_r\qssha
            z_{s_1,\cdots,s_{k-1}}
\\
& =& \sum_{\substack{s_i\geq 1, r,s_1\geq 2\\
    r+s_1+\cdots+s_{k-1}=n}}P(z_{r-1})\qssha
            P(z_{s_1-1,\cdots,s_{k-1}})
\\
    & =& P\big(\hspace{-.3cm}\sum_{\substack{s_i\geq 1, r, s_1\geq 2
                      \\r+s_1+\cdots+s_{k-1}=n
                      }}
                     \hspace{-.3cm} z_{r-1} \qssha P(z_{s_1-1,\cdots,s_{k-1}})\big)
  + P\big( \hspace{-.3cm} \sum_{ \substack{ s_i\geq 1, r, s_1\geq 2 \\
                     r+s_1+\cdots+s_{k-1}=n }}
                     \hspace{-.3cm} P(z_{r-1})\qssha z_{s_1,\cdots,s_{k-1}}\big)
\\
    & = &P(\hspace{-.3cm}\sum_{\substack{ s_i\geq 1, r, s_1\geq 2\\  r+s_1+\cdots+s_{k-1}=n}} \hspace{-.3cm} z_{r-1} \qssha z_{s_1,\cdots,s_{k-1}})
        + P( \hspace{-.3cm} \sum_{\substack{ s_i\geq 1, r, s_1\geq 2\\ r+s_1+\cdots+s_{k-1}=n}} \hspace{-.3cm} z_r\qssha
           z_{s_1-1,\cdots,s_{k-1}}).
\end{eqnarray*}
}
We will denote $\suma$ and $\sumb$ for the first and second sum respectively in the last expression.
Similarly for the sum $\sumz$, we have
{\allowdisplaybreaks
\begin{eqnarray*}
&& \sum_{\substack{r\geq 2, s_i\geq 1\\  r+s_2+\cdots+s_{k-1}=n-1}}z_r\qssha
            z_{1, s_2,\cdots,s_{k-1}}
\\
    &= &\sum_{ \substack{r\geq 2, s_i\geq 1\\
    r+s_2+\cdots+s_k=n-1}}P(z_{r-1})\qssha
            Q(z_{s_2,\cdots,s_{k-1}})
\\
& = &P\bigg(\hspace{-.6cm}\sum_{\substack{r\geq 2, s_i\geq 1\\ r+s_2+\cdots+s_{k-1}=n-1}}
\hspace{-.6cm}z_{r-1}
      \qssha Q(z_{s_2,\cdots,s_{k-1}})\bigg)
     + Q\bigg( \hspace{-.6cm}
     \sum_{\substack{r\geq 2, s_i\geq 1\\ r+s_2+\cdots+s_{k-1}=n-1}}
     \hspace{-.6cm} P(z_{r-1})
       \qssha z_{s_2,\cdots,s_{k-1}} \bigg)
\\
   &=& P(\hspace{-.5cm}
   \sum_{\substack{r\geq 2, s_i\geq 1\\ r+s_2+\cdots+s_{k-1}=n-1}}
   \hspace{-.5cm} z_{r-1} \qssha z_{1,s_2,\cdots,s_{k-1}})
        + Q( \hspace{-.5cm}
        \sum_{\substack{r\geq 2, s_i\geq 1\\ r+s_2+\cdots+s_{k-1}=n-1}}
        \hspace{-.5cm}z_{r} \qssha
           z_{s_2,\cdots,s_{k-1}} )
\end{eqnarray*}
}
We let $\sumc$ and
$\sumd$ to denote the first and the second sum respectively in the last expression.

By Lemma \mref{lem:z1sha} we have
\begin{equation}
\sume= k\sum\limits_{t_i\geq 1, t_1+\cdots+t_{k}=n}
z_{t_1,\cdots,t_{k}} - \sum\limits_{t_i\geq 1,t_k\geq 2,
t_1+\cdots+t_{k}=n} z_{t_1,\cdots,t_{k}}. \mlabel{eq:sum1}
\end{equation}
We also have
{\allowdisplaybreaks
\begin{eqnarray}
\suma+\sumc & =& P(\sum_{\substack{r\geq 2, s_i\geq 1\\ r+s_1+\cdots+s_{k-1}=n}}
z_{r-1}\qssha z_{s_1,s_2,\cdots,s_{k-1}})
\notag \\
& =& P(\sum_{\substack{r, s_i\geq 1\\ r+s_1+\cdots+s_{k-1}=n-1}} z_{r}\qssha
z_{s_1,s_2,\cdots,s_{k-1}})
\notag \\
&= & P\Big(\sum_{\substack{t_i\geq 1\\ t_1+\cdots+t_{k}=n-1}} \calc(t_1,\cdots,
t_{k-1}) z_{t_1,\cdots, t_{k}}\Big)  \ \ \text{(by the induction
hypothesis)}
\label{eq:sum24}
\\
& = &\sum_{ t_i\geq 1, t_1+\cdots+t_{k}=n-1} \calc(t_1,\cdots,
t_{k-1}) z_{t_1+1,\cdots, t_{k}}
\notag\\
& = &\sum_{ t_i\geq 1, t_1\geq 2, t_1+\cdots+t_{k}=n} \calc(t_1-1,
t_2,\cdots, t_{k-1}) z_{t_1,\cdots, t_{k}}.
\notag
\end{eqnarray}
For the sum $\sumb$, we have {\allowdisplaybreaks
\begin{eqnarray}
 \sumb & =& P(\sum_{r\geq 2,s_i\geq 1,
r+s_1+\cdots+s_{k-1}=n-1}z_r\qssha
             z_{s_1,\cdots, s_{k-1}})
\notag \\
     & = &P(\sum_{r , s_i\geq 1, r+s_1+\cdots+s_{k-1}=n-1}z_r\qssha
           z_{s_1,\cdots, s_{k-1}}- \sum_{ s_i\geq 1,
           s_1+\cdots+s_{k-1}=n-2}z_1\qssha z_{s_1,\cdots, s_{k-1}})
\notag
\\
     & =& P\Big(\sum_{t_i\geq 1, t_1+\cdots+t_{k}=n-1}\calc(t_1,\cdots,t_{k-1}) z_{t_1,\cdots,
         t_{k}}-k\sum_{t_i\geq 1, t_k=1, t_1+\cdots+t_{k}=n-1} z_{t_1,\cdots,
         t_{k}}
\notag
\\  &&\quad  -(k-1)\hskip -10 pt\sum_{\substack{t_i\geq 1, t_k\geq 2\\ t_1+\cdots+t_{k}=n-1}}
\hskip -10pt z_{t_1,\cdots,
         t_{k}}\Big) \quad  \text{(by the induction hypothesis and Eq.~(\mref{eq:z1a}))}
\label{eq:sum3}
\\
& =& \sum_{\substack{t_i\geq 1\\
t_1+\cdots+t_{k}=n-1}}\calc(t_1,\cdots,t_{k-1}) z_{t_1+1,\cdots,
         t_{k}}-k\sum_{\substack{t_i\geq 1\\ t_1+\cdots+t_{k}=n-1}} z_{t_1+1,\cdots,t_{k}} \notag
+ \hskip -10pt \sum_{\substack{t_i\geq 1, t_k\geq 2\\
t_1+\cdots+t_{k}=n-1}} \hskip -10pt z_{t_1+1,\cdots,
         t_{k}}
\notag
\\
     & =& \sum_{\substack{
     t_i\geq 1, t_1\geq 2\\t_1+\cdots+ t_{k}=n}}
     \Big(\calc(t_1-1,\cdots, t_{k-1})-k \Big) z_{t_1,\cdots,
     t_{k}} + \sum_{\substack{
     t_i\geq 1, t_1\geq 2, t_k\geq 2\\ t_1+\cdots+ t_{k}=n}}
      z_{t_1,\cdots,
     t_{k}}.
     \notag
\end{eqnarray}
}
For the sum $\sumd$, we have
{\allowdisplaybreaks
\begin{eqnarray}
\sumd &  =& Q\big( \sum_{r, s_i\geq 1, r+s_2+\cdots+s_{k-1}=n-1} z_{r}
            \qssha
           z_{s_2,\cdots,s_{k-1}} -\sum_{s_2+\cdots +s_{k-1}=n-2}z_{1} \qssha
           z_{s_2,\cdots,s_{k-1}}\big)
\notag
\\
      &  =&  Q\Big(\sum_{t_i\geq 1, t_2+\cdots+t_{k}=n-1}
            \calc(t_2, \cdots, t_{k-1})
            z_{t_2,\cdots,t_{k}}
            -\big[\sum_{t_i\geq 1, t_k=1, t_2+\cdots+t_{k}=n-1} (k-1)z_{t_2,\cdots, t_k}
\notag \\
      & & + \hskip -10pt \sum_{\substack{t_i\geq 1, t_k\geq 2 \\ t_2+\cdots+t_{k}=n-1}}
      \hskip -10pt (k-2)z_{t_2,\cdots, t_k} \big]\Big)
      \quad  \text{(by the induction hypothesis and Eq.~(\mref{eq:z1a}))}
\label{eq:sum5}
\\
&  =&  Q\Big(\sum_{t_i\geq 1, t_2+\cdots+t_{k}=n-1}
            [\calc(t_2, \cdots, t_{k-1})-(k-1)]
            z_{t_2,\cdots,t_{k}}+ \hskip -10pt \sum_{\substack{t_i\geq 1, t_k\geq 2 \\
t_2+\cdots+t_{k}=n-1}}
      \hskip -10pt z_{t_2,\cdots, t_k} \Big)
\notag
\\    &  = & \sum_{\substack{ t_i\geq 1, t_1=1 \\t_1+\cdots+t_{k}=n}}
             [ \calc(t_1,\cdots,t_{k-1})- k ] z_{t_1, t_2,\cdots,
              t_{k}}+ \sum_{\substack{ t_i\geq 1, t_1=1, t_k\geq 2
               \\t_1+\cdots+t_{k}=n}}
              z_{t_1, t_2,\cdots,t_{k}} \notag
\end{eqnarray} by Eq.~(\mref{eq:coef2}).

Adding Eq. (\mref{eq:sum24}) and (\mref{eq:sum3}) we obtain
\begin{equation}
\begin{aligned}\suma+\sumb+\sumc & = \hskip -10pt\sum_{\substack{ t_i\geq 1, t_1\geq 2\\
t_1+\cdots+t_{k}=n}}\hskip -10pt (2\calc(t_1\!-\!1,\cdots,
t_{k-1})-k) z_{t_1,\cdots,
     t_{k}} + \hskip -10pt\sum_{\substack{ t_i\geq 1, t_1,t_k\geq 2\\
t_1+\cdots+t_{k}=n}} \hskip -10pt z_{t_1, \cdots, t_k}
     \\
     &= \hskip -10pt \sum_{\substack{t_i\geq 1, t_1\geq 2 \\ t_1+\cdots+
     t_{k}=n}} \hskip -10pt
(\calc(t_1,\cdots, t_{k-1})-k) z_{t_1,\cdots,
     t_{k}}+ \hskip -5pt \sum_{\substack{t_i\geq 1, t_1,t_k\geq 2\\
t_1+\cdots+t_{k}=n}} \hskip -10pt z_{t_1, \cdots, t_k}
     \end{aligned}
     \mlabel{eq:sum234}
\end{equation}
by Eq.~(\mref{eq:coef1}).
Next adding Eq.~(\mref{eq:sum5}) and
(\mref{eq:sum234}) together, we get
\begin{equation}
\suma+\sumb+\sumc+\sumd = \hskip -15pt \sum_{\substack{ t_i\geq 1\\
t_1+\cdots+ t_{k}=n}} \hskip -15pt \Big(\calc(t_1,\cdots,
t_{k-1})-k\Big) z_{t_1,\cdots,
     t_{k}}+
     \hskip -10pt \sum_{\substack{t_i\geq 1, t_k\geq 2\\  t_1+\cdots+ t_{k}=n}}
     \hskip -10pt  z_{t_1,\cdots,
     t_{k}}.
     \mlabel{eq:sum2345}
\end{equation} Finally adding Eq.~(\mref{eq:sum1}) and
(\mref{eq:sum2345}), we obtain Eq.~(\mref{eq:sha}). This completes
the inductive proof of Theorem~\mref{thm:sha}.

\subsection{Proof of Theorem \mref{thm:sha1}.} The proof of Theorem~\mref{thm:sha1} is similar
to that of Theorem \mref{thm:sha}.

First we prove that Eq.~(\mref{eq:sha1}) holds when $k=2$. This follows
from {\allowdisplaybreaks
\begin{eqnarray*} \sum_{\substack{r\geq 1, s_1\geq
2\\r+s_1=n}} z_r\qssha z_{s_1} & =& \sum_{\substack{r\geq 1, s_1\geq
2\\ r+s_1=n}}\ \sum_{\substack{t_1,t_2\geq 1\\
t_1+t_2=n}}\Big(\binc{t_1-1}{r-1}+\binc{t_1-1}{s_1-1}\Big)
z_{t_1}z_{t_2}
\\
 &=& \sum_{t_1,t_2\geq 1, t_1+t_2=n}\bigg (\sum_{r=1}^{\min(t_1,
 n-2)}\binc{t_1-1}{r-1} +\sum_{s_1=2}^{t_1}\binc{t_1-1}{s_1-1}\bigg)z_{t_1}z_{t_2}
 \\
 &= &\Big( \sum_{t_1=1}^{n-2} (2^{t_1}-1)z_{t_1, n-t_1} \Big )+
 (2^{n-1}-2)z_{n-1,1}\\
 &= &\Big( \sum_{t_1=1}^{n-1} (2^{t_1}-1)z_{t_1, n-t_1} \Big )- z_{n-1,1}.
\end{eqnarray*}
}

For the general case we prove Eq.~(\mref{eq:sha1}) by induction on
$n$. If $n=3$, then $k=2$ and we are done. Let $m\geq 4$ be an
integer. Assume that Eq.~(\mref{eq:sha1}) holds when $n\leq m-1$. We
will prove that it holds when $n=m$. Since we have dealt with the
case of $k=2$, we may assume that $k\geq 3$ without loss of
generality.

By Eq.~(\mref{eq:iqform}) the left hand side of Eq.~(\mref{eq:sha1})
is equal to
{\allowdisplaybreaks
\begin{eqnarray*}
&& \sum_{\substack{r, s_i\geq 1, s_1\geq 2\\
r+s_1+\cdots +s_{k-1}=n}}\hskip -20pt  z_{r}\qssha
z_{s_1,\cdots, s_{k-1}}
\\
& =& z_1  \qssha (\sum_{\substack{ s_i\geq 1, s_1\geq 2 \\
s_1+\cdots +s_{k-1}=n-1}}\hskip -20pt z_{s_1,\cdots,
s_{k-1}}) +
\sum_{\substack{r\geq 2, s_i\geq 1, s_1\geq 2 \\
r+s_1+\cdots +s_{k-1}=n}}P(z_{r-1})\qssha
P(z_{s_1-1,s_2,\cdots, s_{k-1}})
\\
&=& z_1  \qssha (\sum_{\substack{ s_i\geq 1, s_1\geq 2 \\
s_1+\cdots +s_{k-1}=n-1}}\hskip -20pt z_{s_1,\cdots,
s_{k-1}}) +
P(\sum_{\substack{r\geq 2, s_i\geq 1, s_1\geq 2 \\
r+s_1+\cdots +s_{k-1}=n}}z_{r-1}\qssha z_{s_1,s_2,\cdots,
s_{k-1}})
\\
&&   + P(\sum_{\substack{r\geq 2, s_i\geq 1, s_1\geq 2 \\
r+s_1+\cdots +s_{k-1}=n}}z_{r}\qssha z_{s_1-1,s_2,\cdots,
s_{k-1}})
\\
&=& z_1  \qssha (\sum_{\substack{ s_i\geq 1, s_1\geq 2 \\
s_1+\cdots +s_{k-1}=n-1}}\hskip -20pt z_{s_1,\cdots,
s_{k-1}}) +
P(\sum_{\substack{r\geq 1, s_i\geq 1, s_1\geq 2 \\
r+s_1+\cdots +s_{k-1}=n-1}}z_{r}\qssha z_{s_1,s_2,\cdots,
s_{k-1}})
\\
&&   + P(\sum_{\substack{r\geq 2, s_i\geq 1 \\
r+s_1+\cdots +s_{k-1}=n-1}}z_{r}\qssha z_{s_1,s_2,\cdots, s_{k-1}}).
\end{eqnarray*}
}

\noindent We denote the three sums in the last expression by
$\sumh$, $\sumi$ and $\sumj$ respectively with the given order.

By Eq.~(\mref{eq:z1b}) and (\mref{eq:coef2}),
we have
{\allowdisplaybreaks
\begin{eqnarray}
\sumh & =& \hskip -10pt \sum_{\substack{ t_i\geq 1,t_1=1, t_2\geq 2\\
t_1+\cdots+t_{k}=n}} \hskip -10pt z_{t_1, t_2,\cdots, t_{k}}+ k
\hskip -10pt \sum_{\substack{t_i\geq 1, t_1\geq 2, t_k=1\\
t_1+\cdots+t_{k}=n}} \hskip -10pt z_{t_1,\cdots, t_{k}}+(k-1)
\hskip -10pt \sum_{\substack{t_i\geq 1, t_1, t_k\geq 2\\
t_1+\cdots+t_{k}=n}}\hskip -10pt z_{t_1,\cdots, t_{k}} \notag
\\
&=& \hskip -10pt \sum_{
\substack{t_i\geq 1, t_1=1\\
t_1+\cdots+t_{k}=n }}\hskip -10pt [ \calc(t_1,\cdots, t_{k})-
\calc(t_2,\cdots, t_{k-1})] z_{t_1,\cdots, t_{k}}- \hskip -10pt
\sum_{ \substack{ t_i\geq 1, t_1=t_2=1
    \\ t_1+\cdots+t_{k}=n }}  z_{t_1,\cdots, t_{k}} \label{eq:sumh}
\\
    &&
+ k
\sum_{ \substack{ t_i\geq 1, t_1\geq 2, t_k=1 \\
t_1+\cdots+t_{k}=n }} z_{t_1,\cdots, t_{k}}+(k-1)
\sum_{ \substack{ t_i\geq 1, t_1\geq 2, t_k\geq 2\\
t_1+\cdots+t_{k}=n }} z_{t_1,\cdots, t_{k}}. \notag
\end{eqnarray}
}
%
By the induction hypothesis, we have
{\allowdisplaybreaks
\begin{eqnarray*}
\sumi &=& P(\sum_{\substack{r\geq 1, s_i\geq 1, s_1\geq 2 \\
r+s_1+\cdots +s_{k-1}=n-1}}z_{r}\qssha z_{s_1,s_2,\cdots,
s_{k-1}})
\\ &= &P\Big(\sum_{
\substack{t_i\geq 1 \\
t_1+\cdots+t_{k}=n-1 }} [\calc(t_1,\cdots,
t_{k-1})-\calc(t_2,\cdots, t_{k-1})] z_{t_1,\cdots, t_{k}} -\sum_{
\substack{ t_i\geq 1, t_2=1
\\
t_1+\cdots+t_{k}=n-1 }}\hskip -20pt   z_{t_1,\cdots,
t_{k}} \Big)
\\
&=& \sum_{
\substack{t_i\geq 1\\
t_1+\cdots+t_{k}=n-1 }} [\calc(t_1,\cdots,
t_{k-1})-\calc(t_2,\cdots, t_{k-1})] z_{t_1+1,\cdots, t_{k}} -
\sum_{ \substack{ t_i\geq 1, t_2=1
\\
t_1+\cdots+t_{k}=n-1 }}\hskip -20pt z_{t_1+1,\cdots,t_{k}}
\\
& =& \sum_{
\substack{  t_i\geq 1,t_1\geq 2 \\
t_1+\cdots+t_{k}=n }} [\calc(t_1-1,\cdots,
t_{k-1})-\calc(t_2,\cdots, t_{k-1})]  z_{t_1,\cdots, t_{k}} - \sum_{
\substack{  t_i\geq 1, t_1\geq 2, t_2=1
\\
t_1+\cdots+t_{k}=n }}  z_{t_1,\cdots, t_{k}}
\end{eqnarray*}
}
By Theorem \mref{thm:sha} and Eq.~(\mref{eq:z1a}), we have
{\allowdisplaybreaks
\begin{eqnarray*}
\sumj &=& P(\sum_{\substack{r\geq 2, s_i\geq 1 \\
  r+s_1+\cdots +s_{k-1}=n-1}}z_{r}\qssha z_{s_1,s_2,\cdots, s_{k-1}})
\\
& =& P(\sum_{\substack{r\geq 1, s_i\geq 1 \\
  r+s_1+\cdots +s_{k-1}=n-1}}z_{r}\qssha z_{s_1,s_2,\cdots,
  s_{k-1}}- \hskip -10pt
  \sum_{\substack{ s_i\geq 1  \\
  s_1+\cdots +s_{k-1}=n-2}} \hskip -10pt z_{1}\qssha z_{s_1,s_2,\cdots, s_{k-1}})
\\
&=& P\Big(\sum_{\substack{t_i\geq 1\\ t_1+\cdots+t_k=n-1}}\calc(t_1,
\cdots, t_{k-1}) -
k\sum_{\substack{t_i\geq 1, t_k=1 \\
t_1+\cdots+t_k=n-1}}z_{t_1,\cdots, t_k}-(k-1)\sum_{\substack{t_i\geq 1, t_k\geq 2 \\
t_1+\cdots+t_k=n-1}}z_{t_1,\cdots, t_k} \Big)
\\
& =& P\Big(\sum_{t_i\geq 1, t_1+\cdots+t_{k}=n-1} [
\calc(t_1,\cdots, t_{k-1}) -k]\: z_{t_1,\cdots,
t_{k}}+\sum_{\substack{t_i\geq 1, t_k\geq 2 \\
t_1+\cdots+t_k=n-1}}z_{t_1,\cdots, t_{k}}\Big)
\\
& =&  \sum_{t_i\geq 1, t_1+\cdots+t_{k}=n-1} [ \calc(t_1,\cdots,
t_{k-1}) -k]\: z_{t_1+1,\cdots, t_{k}}
+\sum_{\substack{t_i\geq 1,t_1,t_k\geq 2\\
t_1+\cdots+t_k=n-1 }}z_{t_1+1,\cdots, t_{k}}
\\
&=& \sum_{t_i\geq 1, t_1\geq 2,  t_1+\cdots+t_{k}=n}
[\calc(t_1-1,t_2,\cdots, t_{k-1}) -k] \: z_{t_1,\cdots, t_{k}}
+\sum_{\substack{t_i\geq 1,t_1,t_k\geq 2\\
t_1+\cdots+t_k=n}}z_{t_1,\cdots, t_{k}}.
\end{eqnarray*}
}
Hence, using
Eq.~(\mref{eq:coef1}) we obtain
{ \allowdisplaybreaks
\begin{eqnarray}
\sumi+\sumj & =& \sum_{\substack{ t_i\geq 1, t_1\geq 2
\\ t_1+\cdots+t_{k}=n
}} \hskip -10pt \Big(2\calc(t_1-1,t_2,\cdots,
t_{k-1})-\calc(t_2,\cdots, t_{k-1})-k \Big)z_{t_1,\cdots, t_{k}}
\notag \\
& & - \sum_{  \substack{ t_i\geq 1, t_1\geq 2, t_2=1
\\
t_1+\cdots+t_{k}=n
}
} \hskip -20pt z_{t_1,\cdots, t_{k}}+\sum_{\substack{t_i\geq 1,t_1,t_k\geq 2\\
t_1+\cdots+t_k=n}}z_{t_1,\cdots, t_{k}} \label{eq:sum78}
\\
&=& \sum_{\substack{ t_i\geq 1, t_1\geq 2 \\
t_1+\cdots+t_{k}=n}} \Big(\calc(t_1,t_2,\cdots,
t_k)-\calc(t_2,\cdots, t_{k-1})-k \Big)z_{t_1,\cdots, t_{k}}
\notag \\
& &  -\sum_{\substack{t_i\geq 1, t_1\geq 2, t_2=1 \\
t_1+\cdots+t_{k}=n }} z_{t_1,\cdots, t_{k}} +\sum_{\substack{t_i\geq 1,t_1,t_k\geq 2\\
t_1+\cdots+t_k=n}}z_{t_1,\cdots, t_{k}}. \notag
\end{eqnarray}
}
Adding Eq.~(\mref{eq:sumh}) and (\mref{eq:sum78}) together we obtain
\begin{eqnarray*} \sumh+\sumi+\sumj
 &=&   \sum_{
\substack{t_i\geq 1, t_1=1 \\
t_1+\cdots+t_{k}=n }} [ \calc(t_1,\cdots,
t_{k})- \calc(t_2,\cdots, t_{k-1})] z_{t_1,\cdots, t_{k}} - \hskip
-10pt \sum_{ \substack{ t_i\geq 1, t_1=t_2=1
    \\
        t_1+\cdots+t_{k}=n }} \hskip -10pt z_{t_1,\cdots, t_{k}}
\\ & &+ k
\sum_{ \substack{ t_i\geq 1, t_1\geq 2, t_k=1 \\
t_1+\cdots+t_{k}=n }} z_{t_1,\cdots, t_{k}}+ (k-1)
\sum_{ \substack{ t_i\geq 1, t_1\geq 2,t_k\geq 2 \\
t_1+\cdots+t_{k}=n }} z_{t_1,\cdots, t_{k}}
\\ &&+ \sum_{\substack{ t_i\geq 1,t_1\geq 2 \\
t_1+\cdots+t_{k}=n}} [\calc(t_1,t_2,\cdots,
t_k)-\calc(t_2,\cdots, t_{k-1})-k]z_{t_1,\cdots, t_{k}}
\\
 &&
-\sum_{\substack{t_i\geq 1, t_1\geq 2, t_2=1 \\
t_1+\cdots+t_{k}=n }} z_{t_1,\cdots, t_{k}} +\sum_{\substack{t_i\geq 1,t_1,t_k\geq 2\\
t_1+\cdots+t_k=n}}z_{t_1,\cdots, t_{k}}
\\
&=&  \hskip -10pt \sum_{
\substack{t_i\geq 1, t_1=1 \\
t_1+\cdots+t_{k}=n }}\hskip -10pt [ \calc(t_1,\cdots,
t_{k})- \calc(t_2,\cdots, t_{k-1})] z_{t_1,\cdots, t_{k}}
 - \hskip -10pt \sum_{ \substack{ t_i\geq 1, t_1=t_2=1
    \\
        t_1+\cdots+t_{k}=n }}\hskip -10pt z_{t_1,\cdots, t_{k}}
\\
&&+ \sum_{\substack{ t_i\geq 1,t_1\geq 2 \\
t_1+\cdots+t_{k}=n}} [\calc(t_1,t_2,\cdots,
t_k)-\calc(t_2,\cdots, t_{k-1})]z_{t_1,\cdots, t_{k}}
 -\hskip -10pt \sum_{\substack{  t_i\geq 1,t_1\geq 2, t_2=1\\
t_1+\cdots+t_{k}=n }} \hskip -10pt z_{t_1,\cdots, t_{k}}
\\
&=& \sum_{\substack{ t_i\geq 1 \\
t_1+\cdots+t_{k}=n}} [\calc(t_1,t_2,\cdots,
t_k)-\calc(t_2,\cdots, t_{k-1})]z_{t_1,\cdots, t_{k}}-
\hskip -10pt  \sum_{\substack{  t_i\geq 1, t_2=1\\
t_1+\cdots+t_{k}=n }} \hskip -10pt z_{t_1,\cdots, t_{k}},
\end{eqnarray*}
as desired. This completes the inductive proof of Theorem~\mref{thm:sha1}.


\end{document}